\documentclass[11pt]{amsart}

\usepackage[square,compress,comma, numbers,sort]{natbib}
\usepackage[colorlinks=true, citecolor=red, linkcolor=blue]{hyperref}
\usepackage{amsfonts,mathtools}

\usepackage{caption}
\usepackage{subcaption}

\allowdisplaybreaks[4]

\usepackage[T2A]{fontenc}
\usepackage[utf8]{inputenc}
\usepackage{polski}
\usepackage[english]{babel}

\usepackage{graphicx}
\usepackage{float}
\graphicspath{}

\DeclareGraphicsExtensions{.pdf,.png,.jpg}
\usepackage{amsmath}
\usepackage{bbm}
\usepackage{amssymb}
\usepackage{mathtools}

\usepackage{color}
\definecolor{c20}{rgb}{0.,0.7,0.}
\definecolor{c30}{rgb}{0.,0.,1.}
\definecolor{c40}{rgb}{1,0.1,0.7}
\definecolor{c50}{rgb}{1,0,0}
\definecolor{c60}{rgb}{1,0.9,0.1}

\def\Var{\text{Var}}

\definecolor{c20}{rgb}{0.,0.7,0.}
\definecolor{c30}{rgb}{0.,0.,1.}
\definecolor{c40}{rgb}{1,0.1,0.7}
\definecolor{c50}{rgb}{1,0,0}
\definecolor{c60}{rgb}{1,0.9,0.1}

\newcommand{\ve}{\varepsilon}

\newcommand{\E}[1]{\mathbb{E}\left\{ #1\right\}}

\newcommand{\pk}[1]{\mathbb{P} \left\{ #1 \right \} }

\newcommand{\R}{\mathbb{R}}

\newcommand{\inr}{\in \R}

\newcommand{\mD}{\mathbb{D}}

\newcommand{\limit}[1]{\lim_{#1 \to   \infty}}

\newcommand{\BQN}{\begin{eqnarray}}
\newcommand{\EQN}{\end{eqnarray}}
\newcommand{\BQNY}{\begin{eqnarray*}}
\newcommand{\EQNY}{\end{eqnarray*}}

\newcommand{\BS}{\begin{sat}}
\newcommand{\ES}{\end{sat}}
\newcommand{\BT}{\begin{theo}}
\newcommand{\ET}{\end{theo}}
\newcommand{\BK}{\begin{korr}}
\newcommand{\EK}{\end{korr}}

\newcommand{\PH}{\overline{\Phi}}
\DeclareMathOperator{\cov}{cov}

\newcommand{\BD}{\begin{de}}
\newcommand{\ED}{\end{de}}

\newcommand{\BIT}{\begin{itemize}}
\newcommand{\EIT}{\end{itemize}}
\newcommand{\BDI}{\begin{description}}
\newcommand{\EDI}{\end{description}}

\newcommand{\BRM}{\begin{remarks}}
\newcommand{\ERM}{\end{remarks}}

\newcommand{\BEL}{\begin{lem}}
\newcommand{\EEL}{\end{lem}}

\def\bqny#1{{\begin{eqnarray*} #1 \end{eqnarray*}}}
\def\bqn#1{{\begin{eqnarray} #1 \end{eqnarray}}}

\newtheorem{theo}{Theorem}[section]
\newtheorem{sat}[theo]{Proposition}
\newtheorem{de}[theo]{Definition}
\newtheorem{lem}[theo]{Lemma}

\newtheorem{korr}[theo]{Corollary}

\newtheorem{remarks}[theo]{Remarks}
\newtheorem{prop}[theo]{Proposition}

\newcommand{\COM}[1]{}

\def\IF{\infty}

\newcommand{\QED}{\hfill $\Box$}

\def\ve{\varepsilon}

\def\IF{\infty}

\newcommand{\pr}{\mathcal{P}_{T_u}}
\newcommand{\prt}{\mathcal{P}_{T}}

\newcommand{\ind}{\mathbb{I}}

\begin{document}

\title{ Parisian Ruin for Insurer and Reinsurer under Quota-Share
Treaty}

\author{{Grigori Jasnovidov}}
\address{Grigori Jasnovidov,
St. Petersburg Department
of Steklov Mathematical Institute
of Russian Academy of Sciences, St. Petersburg, Russia
}
\email{griga1995@yandex.ru}

\author{{Aleksandr Shemendyuk}}
\address{Aleksandr Shemendyuk, Department of Actuarial Science, 
University of Lausanne,\\
UNIL-Dorigny, 1015 Lausanne, Switzerland
}
\email{Aleksandr.Shemendyuk@unil.ch}

\bigskip

\date{\today}
 \maketitle

{\bf Abstract:}
In this contribution we study asymptotics of the simultaneous
Parisian ruin probability of a
two-dimensional fractional Brownian motion risk process.
This risk process models the surplus processes of
an insurance and a reinsurance companies, where the net loss is distributed between them in given proportions.
We also propose an approach for simulation of Pickands and Piterbarg
type constants appearing in the asymptotics of the ruin probability.



\section{Introduction}

Consider the risk model defined by
\bqn{\label{Risk_Model}
R(t) = u+\rho t-X(t), \ \ \ \ t \ge 0,
}
where $X(t)$ is a centered Gaussian risk process
with a.s. continuous sample paths,
$\rho>0$ is the net profit rate and $u>0$ is the initial capital.
This model is relevant to insurance and financial applications, see,
e.g, \cite{MR1458613}.
A question of numerous investigations
is study of the asymptotics of the classical ruin
probability
\bqn{\label{classical_ruin_probability}
\lambda(u) := \pk{\exists t\ge 0: R(t)<0}
}
as $u \to \IF$ under different levels of generality. It turns out, that only
for $X$ being a Brownian motion (later on BM) $\lambda(u)$
can be calculated explicitly. Namely,
if $X$ is a standard BM, then $\lambda(u) = e ^{-2\rho u}, \ u,\rho>0,$ see,
e.g., \cite{DeM15}. Since it seems impossible to find the exact value
of $\lambda(u)$ in other cases, asymptotics of $\lambda(u)$ as $u \to \IF$
is dealt with. First the
problem of a large excursion of a stationary Gaussian process
was considered by J. Pickands in 1969, see \cite{PickandsA}.
We refer to monographs \cite{20lectures,Pit96,MR1379077}
for the survey of known results by the recent time.
We would like to point out seminal manuscript \cite{DI2005}
establishing asymptotics of $\lambda(u)$
under week assumptions on variance and covariance of $X$.
For the
discrete-time investigations (i.e., when $t$ in model \eqref{Risk_Model}
belongs to a discrete grid $\{0,\delta,2\delta...\}$ for some $\delta>0$),
we refer to \cite{HashorvaGrigori,Piterbargdiscrete,
secondproj,grisha_kdebicki_gaussian_reflected,pick_speed_convergence,
piterbarg_speed_convergence}.
We would like to suggest a reader contributions \cite{KrzysPeng2015,
DebickiScaledBM,HashDebPenggammareflect,SojournInfty,
ParisianFiniteHoryzon,SojournHashorvaCorrelatedBM,SumFBMDebicki,
Rolski17,Bai_2017,ParisRuinGenrealHashorvaJiDebicki}
for the related generalizations of the classical ruin problem.
Some contributions (see, e.g., \cite{ParisianFiniteHoryzon,
ParisRuinGenrealHashorvaJiDebicki,Bai_2017}), extend the
classical ruin problem to the so-called Parisian ruin problem which
allows the surplus process to spend a pre-specified time below zero before a ruin is recognized.
Formally, the classical Parisian ruin probability is defined by
\bqn{\label{paris_ruin-def}
\pk{\exists t\ge 0 : \forall s\in [t,t+T] \ R(s)<0}, \ \ \ \ \ \ \ \
T\ge 0.
}
As in the classical case, only for $X$ being a BM
 the probability above can
be calculated explicitly (see \cite{LCPalowski}):
$$\pk{\exists t\ge 0 : \forall s\in [t,t+T] \ B(s)-cs>u} =
\frac{e^{-c^2T/2}-c\sqrt {2\pi T}\Phi(-c\sqrt T)}{
e^{-c^2T/2}+c\sqrt {2\pi T}\Phi(c\sqrt T)}e^{-2cu}
, \ \ \ T\ge 0$$
where $\Phi$ is the distribution function of a standard Gaussian random variable and $B$ is a standard BM.
Note in passing, that the asymptotics of the Parisian ruin probability
for $X$ being a self-similar Gaussian processes is derived in
\cite{ParisRuinGenrealHashorvaJiDebicki}.
We refer to \cite{HashorvaGrigori,ParisianFiniteHoryzon} for investigations of some other problems in this field.
\\

Motivated by
\cite{Lanpeng2BM} (see also \cite{secondproj,sojourn_2_dim}),
we study a model where two companies share the net losses in proportions $\delta_1,\delta_2 > 0$, with $\delta_1 +
\delta_2 = 1$, and receive the premiums at rates $\rho_1 , \rho_2 > 0$, respectively. Further, the risk process of the $i$th company is defined by
$$R_i (t) = x _i + \rho_i t - \delta _i B (t), \ \ \  t \ge 0,  \ i = 1,2,$$
where $x_i > 0$ is the initial capital of the $i$th company. In this model both claims and net losses are distributed between the companies, which corresponds to the proportional reinsurance dependence of the companies.
In this paper we study the asymptotics of
the simultaneous Parisian ruin probability defined by
$$\pk{\exists t\ge 0: \forall s \in [t,t+T] \ R_1(s)<0,R_2(s)<0}
, \ \ \ T\ge 0.$$
Since the probability above
does not change under a scaling of $(R_1 ,R_2)$,
it equals to
$$\pk{\exists t\ge 0: \forall s \in [t,t+T] \ u_1+c_1s-B(s)<0,
u_2+c_2s-B(s)<0}, \ \ \ T\ge 0,$$
where $u_i = x_i/\delta_i$ and $ c_i = \rho_i/\delta_i, \ i = 1,2$.
Later on we derive the asymptotics of
the probability above as $u_1 ,u_2$ tend to infinity at the constant speed (i.e., $u_1 /u_2$ is constant). Therefore, we let $u _i = q _i u$ be fixed constants with $q _i > 0, \ i = 1,2$
and deal with asymptotics of
$$\mathcal P_T(u):=
\pk{\exists t\ge 0: \forall s \in [t,t+T] \   B(s)>q_1u+c_1s,B(s)>q_2u+c_1s}
, \ \ \ T\ge 0$$
as $u\to \IF$.
Letting the initial capital tends to infinity is not just a mathematical assumption, but also an economic requirement stated by authorities in all developed countries,
see \cite{Mikosch2004}.
It aims to prevents a company from the bankruptcy
because of excessive number of small claims and/or several major claims,
before the premium income is able to balance the losses and profits.
Observe that $\prt(u)$ can be rewritten as
$$\pk{\exists t\ge 0:\forall s\in [t,t+T]\  B(s)-\max(c_1s+q_1u,c_2s+q_2u)>0
}.$$

Thus, the two-dimensional problem may also be considered
as a one-dimensional crossing problem over a piece-wise linear barrier.
If the two lines $q _1 u + c _1 t$ and $ q_ 2 u + c _2 t$ do
not intersect over $(0,\IF)$, then the problem reduces
to the classical one-dimensional BM risk model, which has been
discussed in \cite{ParisRuinGenrealHashorvaJiDebicki,ParisianFiniteHoryzon}
 and thus will not be the focus of this paper.
In consideration of that, we assume that
\bqn{\label{cq}
c_1>c_2, \ \ \ q_2>q_1.
}
Under the assumption above the
lines $q _1 u + c _1 t$ and $q_ 2 u + c _2 t$ intersects
at point $ut_*$ with
\bqn{\label{t_*_new}
t_* = \frac{q_2-q_1}{c_1-c_2}>0
}
that plays a crucial role in the following.
The first usual step when dealing with asymptotics of a ruin
probability of a Gaussian process is centralizing the process. In
our case it can be achieved by the self-similarity of BM:
\bqny{
\prt(u) &=&
\pk{\exists tu\ge 0:\inf\limits_{su\in [tu,tu+T]}(B(su)-c_1su)>q_1u,
\inf\limits_{su\in [tu,tu+T]}(B(su)-c_2su)>q_2u}
\\&=&
\pk{\exists t\ge 0:\inf\limits_{s\in [t,t+T/u]}
(B(s)-(c_1s+q_1) \sqrt u )>0,
\inf\limits_{s\in [t,t+T/u]}(B(s)-(c_2s+q_2) \sqrt u)>0}
\\&=&
\pk{\exists t\ge 0:\inf\limits_{s\in [t,t+T/u]}\frac{B(s)}
{\max(c_1s+q_1,c_2s+q_2)}> \sqrt u}
.}
The next step is analysis of the variance of the centered process.
Note that the variance of
$\frac{B(t)}{\max(c_1t+q_1,c_2t+q_2)}$ can achieve its unique maxima only
at one of the following points:
\bqny{\label{ttdefinition_new}
t_*, \ \ \ \overline t_1 := \frac{q_1}{c_1}, \ \ \
\overline t_2 := \frac{q_2}{c_2}.
}
From \eqref{cq} it follows that $\overline t_1<\overline t_2$.
As we see later, the position of $t_*$ regarding to $(\overline t_1,
\overline t_2)$ determines the asymptotics of $\prt(u)$.
Note, that the variance of
$\frac{B(t)}{\max(c_1t+q_1,c_2t+q_2)}$
 is not smooth around  $t_*$
 if \eqref{cq} is satisfied; this observation does not allow us to obtain the asymptotics
of $\prt(u)$ straightforwardly by \cite{ParisRuinGenrealHashorvaJiDebicki}.
Define for any $L\ge 0$ and some function $h: \R \to \R$ constant
\bqny{
\mathcal{F}_L^{h} = \E{\sup\limits_{t\in \R}\inf\limits_{s\in [t,t+L]}
e^{\sqrt 2 B(s)-|s|+h(s)}}
}
when the expectation above is finite. For the properties of
 $\mathcal{F}_L^{h}$ and related Piterbarg constants
 we refer to
\cite{ParisianFiniteHoryzon, ParisRuinGenrealHashorvaJiDebicki,
ApproximationSupremumMaxstableProcess,20lectures}.
Notice that $\mathcal{F}_0^{h}$ is the Piterbarg
constant introduced in \cite{Lanpeng2BM}.
Let $\PH$ be the survival function of a standard Gaussian random
variable and $\mathbb{I}({\cdot})$ be the indicator function.
The next theorem derives the asymptotics of $\prt(u)$ as
$u \to \IF$:
\begin{theo} \label{theoparis_simple}
 Assume that \eqref{cq} holds.\\
1)If $t_*\notin(\overline t_1,\overline t_2)$, then as $u\to \IF$
\bqn{\label{parisclaim11}
\prt(u) \sim
\left(\frac{1}{2}\right)^{\mathbb{I}(t_*=\overline t_i)}
\frac{e^{-c_i^2T/2}-c_i\sqrt{2\pi T}\Phi(-c_i\sqrt T)}
{e^{-c_i^2T/2}+c_i\sqrt{2\pi T}\Phi(c_i\sqrt T)}e^{-2c_iq_iu},
}
where $i=1$ if $t_*\le \overline t_1$ and $i=2$ if
$t_*\ge \overline t_2$.

2)If $t_*\in (\overline t_1,\overline t_2)$, then
 as $u\to \IF$
\bqny{
\prt(u) \sim
\mathcal{F}_{T'}^d\PH\left(
(c_1q_2-c_2q_1)\sqrt \frac{q_2-q_1}{c_1-c_2} \sqrt u\right),
}
where $\mathcal{F}_{T'}^d \in (0,\IF)$ and
\bqn{\label{defT'd}
\ \ \ \ T' = T\frac{(c_1q_2-q_1c_2)^2}{2(c_1-c_2)^2}, \ \ \
\ \ \ \ d(s) =
s\frac{c_1q_2+c_2q_1-2c_2q_2}{c_1q_2-q_1c_2}\mathbb{I}(s<0)
+s\frac{2c_1q_1-c_1q_2-q_1c_2}{c_1q_2-q_1c_2}\mathbb{I}(s\ge0)
.}
\end{theo}

\COM{\\
\\
We organize the rest of the paper in the following way.
The next section gives all necessary notation and the main results.
All proofs are relegated to Section 3 followed by an Appendix containing
auxiliary calculations.}

\section{Main Results}

In classical risk theory, the surplus process of an insurance company is modeled by the compound Poisson or the general compound renewal risk process, see, e.g., \cite{MR1458613}. The calculation of the ruin probabilities is of a particular interest for both theoretical and applied domains. To avoid the technical issues and allow for dependence between claim sizes, these models are often approximated by the risk model \eqref{Risk_Model},  driven by $B_H$ a
standard fractional Brownian motion (later on fBm), i.e, Gaussian process with zero-mean and covariance function
$$\cov(B_H(t),B_H(s))=\frac{t^{2H}+s^{2H}-|t-s|^{2H}}{2}, \ \ \ \ s,t \in
\R, \ H \in (0,1).$$
Since the time spent by the surplus process below zero may depend on $u$, in the following we allow $T=:T_u$ in \eqref{paris_ruin-def}
to depend on $u$. As mentioned in
\cite{ParisianFiniteHoryzon},
for the one-dimensional Parisian ruin probability we need to control
the growth of $T_u$ as $u \to \IF$. Namely, we impose the following condition:
\bqn{\label{assumption_T_u}
 \limit{u}T_u u^{1/H-2} =T \in [0,\IF), \ H\in (0,1)
.}
Note that if $H>1/2$, then $T_u$ may grow to infinity, while
if $H<1/2$, then $T_u$ approaches zero as $u$ tends to infinity.
As we see later in Proposition \ref{remark_H<1/2}, the condition above is
necessary and the result does not hold without it.
As for BM, by the self-similarity of fBm we obtain
\bqny{
\pr(u) &:=&
\pk{\exists t\ge 0: \forall s \in [t,t+T_u] \   B_H(s)>q_1u+c_1s,
B_H(s)>q_2u+c_1s}
\\&=&
\pk{\exists t\ge 0:\inf\limits_{s\in [t,t+T_u/u]}\frac{B_H(s)}
{\max(c_1s+q_1,c_2s+q_2)}> u^{1-H}}
.}
The variance of
$\frac{B_H(t)}{\max(c_1t+q_1,c_2t+q_2)}$ can achieve its unique maxima only
at one of the following points:
\bqn{\label{ttdefinition}
t_*, \ \ \ t_1 := \frac{Hq_1}{(1-H)c_1}, \ \ \ t_2 := \frac{Hq_2}{(1-H)c_2}.
}

From \eqref{cq} it follows that $t_1<t_2$.
Again, the position of $t_*$ regarding to $(t_1,t_2)$ determines
the asymptotics of $\pr(u)$.
Define for $H\in (0,1)$ and $T\ge 0$ Pickands constants by
$$\mathbb{H}_{2H} =\limit{S} \frac{1}{S}\E{\sup\limits_{t \in[0,S]}
e^{\sqrt 2 B_H(t)-t^{2H}}}, \ \ \
\ \
\mathcal{F}_{2H}(T) = \limit{S}\frac{1}{S}
\E{\sup\limits_{t\in [0,S]}\inf\limits_{s\in [0,T]}
e^{\sqrt 2 B_H(t+s)-(t+s)^{2H}}
}.$$
It is shown in \cite{ParisRuinGenrealHashorvaJiDebicki} and
\cite{20lectures}, respectively, that
$\mathcal{F}_{2H}(T)$ and $\mathbb{H}_{2H}$ are finite positive
constants.
Let
\bqn{\label{mathbb_D}
 \ \ \mD_H = \frac{c_1t_*+q_1}{t_*^H}, \
K_H = \frac{2^{\frac{1}{2}-\frac{1}{2H}}\sqrt \pi}{\sqrt{H(1-H)}},
\
\mathbb{C}_{H}^{(i)} = \frac{c_i^Hq_i^{1-H}}{H^H(1-H)^{1-H}},
\ D_i =
\frac{c_i^2(1-H)^{2-\frac{1}{H}}}{2^{\frac{1}{2H}}H^{2}}, \ i=1,2.
}

Now we are ready to give the asymptotics of $\pr(u)$:
\begin{theo} \label{theoparis}
 Assume that \eqref{cq} holds and $T_u$ satisfies \eqref{assumption_T_u}.\\
1)If $t_*\notin(t_1,t_2)$, then as $u\to \IF$
\bqn{\label{parisclaim11}
\pr(u) \sim
\left(\frac{1}{2}\right)^{\mathbb{I}(t_*=t_i)}\times
\begin{cases}
\frac{e^{-c_i^2T/2}-c_i\sqrt{2\pi T}\Phi(-c_i\sqrt T)}
{e^{-c_i^2T/2}+c_i\sqrt{2\pi T}\Phi(c_i\sqrt T)}e^{-2c_iq_iu},
&H=1/2,\\
K_H\mathcal{F}_{2H}(T D_i)(\mathbb{C}_H^{(i)}u^{1-H})^{\frac{1}{H}-1}
\PH(\mathbb{C}_H^{(i)}u^{1-H}),
& H\neq 1/2,
\end{cases}}
where $i=1$ if $t_*\le t_1$ and $i=2$ if $t_*\ge t_2$.

2)If $t_*\in (t_1,t_2)$ and $\lim\limits_{u\to \IF}T_uu^{2-1/H} = 0$ for
$H>1/2$, then as $u \to \IF$
\bqn{\label{parisclaim2}
\pr(u) \sim
\PH(\mathbb{D}_{H} u^{1-H})\times
\begin{cases}
1, &H>1/2,\\
\mathcal{F}_{T'}^d,
&H=1/2,\\
\mathcal{F}_{2H}(\overline D T)
 Au^{(1-H)(1/H-2)},& H<1/2,
\end{cases}}
where $\mathcal{F}_{T'}^d \in (0,\IF)$ with $T'$ and $d$
 defined in \eqref{defT'd} and
\bqn{\label{A_overline_D} \ \
A = \Big(|H(c_1t_*+q_1)-c_1t_*|^{-1} +
|H(c_2t_*+q_2)-c_2t_*|^{-1}\Big) \frac{t_*^H\mD_H^{\frac{1}{H}-1}}
{2^{\frac{1}{2H}}}, \ \ \ \ \ \
\overline D=\frac{(c_1t_*+q_1)^{\frac{1}{H}}}{2^{
\frac{1}{2H}
}t_*^2}. }
\end{theo}
The theorem above generalizes Theorem \ref{theoparis_simple}
and Theorem 3.1 in \cite{Lanpeng2BM}. Note that if $T=0$, then
the result above reduces to Theorem 3.1 in \cite{Lanpeng2BM}.
As indicated in \cite{ParisianFiniteHoryzon}, it seems extremely difficult to find the exact asymptotics of the
one-dimensional Parisian ruin probability
 if \eqref{assumption_T_u} does not hold.
The intuitive reason is that the ruin happens over
'too long interval'. To illustrate difficulties arising in approximation
of $\pr(u)$ in this setup we consider a 'simple'
scenario: let $T_u=T>0$ and $H<1/2$. In this case we have
\begin{prop} \label{remark_H<1/2}
If $H<1/2, \ T_u=T>0$ and $t_* \in (t_1,t_2)$, then
\bqn{\label{paris_ineq_Hsmall}
\ \ \ \
\bar{C} \PH(\mD_H u^{1-H})e^{-C_{1,\alpha} u^{2-4H}-
 C_{2,\alpha}u^{2(1-3H)}}
\le \pr(u)
\le
(2+o(1))\PH(\mD_Hu^{1-H})\PH\left(u^{1-2H}\frac{T^H\mD_H}{2t_*^H}\right)
,}
where $\bar C\in (0,1)$ is a fixed constant that does not depend on $u$ and
\bqn{\label{C_alpha_1,2}
\alpha = \frac{T^{2H}}{2t_*^{2H}}, \ \ \ \ \ \
C_{i,\alpha} = \frac{\alpha^{i}}{i} \mD_H^2, \ \ \
i=1,2
.}
\end{prop}
Note that the result above expands Theorem 3.2
in \cite{ParisianFiniteHoryzon} for fBm case. Comparing 
the proposition above with Theorem \ref{theoparis} we see that right hand part above is exponentially smaller
than the corresponding asymptotics in case $H<1/2$ in
Theorem \ref{theoparis} and hence condition
\eqref{assumption_T_u} indeed is important.

\section{Simulation of Piterbarg \& Pickands constants}
In this section we give algorithms for simulations of Pickands and Piterbarg
type constants appearing in Theorems
\ref{theoparis_simple} and \ref{theoparis} and study their properties
relevant for simulations.
Since the classical Pickands constant $\mathbb H_{2H}$ has been investigated
in several contributions (see, e.g., \cite{DiekerY} and references therein),
later on we deal with $\mathcal F_L^h$ and $\mathcal F_{2H}(L)$.
For notation simplicity we denote for any real numbers $x<y$ and $\tau>0$
$$[x,y]_\tau = [x,y]\cap \tau\mathbb Z.$$

\textbf{Simulation of Piterbarg constant.}
In this subsection we always assume that
$$L\ge 0 \ \ \ \text{and} \ \ \
h(s) = b s \  \ind (s<0)-a s \ \ind(s\ge 0), \ \ s\inr, \ a,b>0.$$
To simulate $\mathcal{F}_L^{h}$ we use approximation
\bqny{
\mathcal{F}_L^{h} \approx
\E{\sup\limits_{t\in [-M,M]_{\tau}}\inf\limits_{s\in [t,t+L]_{\tau}}
e^{\sqrt 2 B(s)-|s|+h(s)}},
}
where $M$ is sufficiently large and $\tau>0$ is sufficiently small.
The approximation above has several errors: truncation error
(i.e, choice of $M$), discretization error (i.e., choice of $\tau$)
and simulation error.
It seems difficult to give a precise estimate of the discretization
error, we refer to \cite{DiekerY,pick_speed_convergence} for discussion of such problems.
To take an appropriate $M$
and give an upper bound of the truncation error we
derive few lemmas. The first lemma provides us bounds for
$\mathcal F_L^h$:
\begin{lem} \label{lower_bound_for_Piterb_const}
It holds that for $L>0$
\bqny{
 2e^{- L\min(a,b)}\PH(\sqrt{2L})\le \mathcal{F}_L^{h}
<1+\frac{1}{a}+\frac{1}{b}-\frac{1}{a+b+1}
}
and
\bqny{
\mathcal{F}_0^{h} = 1+\frac{1}{a}+\frac{1}{b}-\frac{1}{a+b+1}.
}
\end{lem}
Note that the second statement of the lemma gives
the explicit expression for the two-sided Piterbarg constant introduced
in \cite{Lanpeng2BM}.
In the next lemma we focus on the truncation error:
\begin{lem}\label{lemma_estimation_of_supinf_over_infty}
For $M\ge 0$ it holds that
\bqn{\label{estimation_of_supinf_over_infty}
\E{\sup\limits_{t\in \R \backslash [-M,M]}\inf\limits_{s\in [t,t+L]}
e^{\sqrt 2 B(s)-|s|+h(s)}} \le e^{-aM}\left(1+\frac{1}{a}\right)+
e^{-bM}\left(1+\frac{1}{b}\right).
}
\end{lem}

Now we are ready to find an appropriate $M$. We have
by Lemma \ref{lemma_estimation_of_supinf_over_infty} that
\bqny{
\left|\mathcal{F}_L^{h} -
\E{\sup\limits_{t\in [-M,M]}\inf\limits_{s\in [t,t+L]}
e^{\sqrt 2 B(s)-|s|+h(s)}}\right| &\le&
\E{\sup\limits_{t\in \R \backslash [-M,M]}\inf\limits_{s\in [t,t+L]}
e^{\sqrt 2 B(s)-|s|+h(s)}}
\\ &\le&
2\left(1+\frac{1}{\min(a,b)}\right)e^{-M\min(a,b)}
}
and on the other hand by Lemma \ref{lower_bound_for_Piterb_const}
$$\mathcal{F}_L^{h}\ge 2e^{- L\min(a,b)}\PH(\sqrt{2L}),$$
hence to obtain a good accuracy we need that
$$\left(1+\frac{1}{\min(a,b)}\right)e^{-\min(a,b)M} \ll e^{-L\min(a,b)}\PH(\sqrt{2L}).$$
Assume for simulations that $\min(a,b)\ge 1$; otherwise special
case $\min(a,b)<<1$ requires a choice of a large $M$ implying
very high level of computation capacity.
For simulations, we take $M  = \frac{7+L(3+\min(a,b))}{\min(a,b)}$
providing us truncation error smaller than $3*10^{-3}$;
we do not need to have better due to the discretization and simulation
errors.
Since it seems difficult to estimate these errors,
 we just take a 'small' $\tau$ and  a 'big' number of simulation $n$.
The above observations give us the following algorithm:
\\\\
1) take $M = \frac{7+L(3+\min(a,b))}{\min(a,b)},
\tau = 0.005$ and $n=10^4$;\\
2) simulate $n$ times $B(t), \ t \in [-M,M]_\tau$,
i.e, obtain $B_i(t), \ 1\le i\le n$;\\
3) compute
 \bqny{
 \widehat{\mathcal F_L^h}:=  \frac{1}{n}\sum\limits_{i=1}^{n}
\sup\limits_{t\in [-M,M]_\tau}\inf\limits_{s\in [t,t+L]_\tau}
e^{\sqrt 2 B_i(s)-|s|+h(s)}.
}

\textbf{Simulation of Pickands constant.}
It seems difficult to simulate $\mathcal{F}_{2H}(L)$ relying
straightforwardly on its definition. As follows from approach in \cite{SBK,DiekerY} for any $\eta>0$ with $W(t) = B_{2H}(t)-|t|^{2H}$
$$\mathcal F_{2H}(L)
=
\E{ \frac{\sup\limits_{t\in \R}
\inf\limits_{s\in [t,t+L]} e^{W(t)}}
{\eta \sum\limits_{k\in \mathbb Z} e^{W(k\eta)}}}
 .$$
The merit of the representation above is that there is no limit as
is in the original definition and thus it is much easier to simulate
$\mathcal F_{2H}(L)$ by the Monte-Carlo method. The second benefit is that
there is a sum in the denominator, that can be
simulated easily with a good accuracy. The only drawback is that
 the $\sup\inf$ in the nominator is taken on the whole real line.
Thus, we approximate $\mathcal F_{2H}(L)$ by discrete analog of
the formula above:
$$\mathcal F_{2H}(L) \approx
\E{ \frac{\sup\limits_{t\in [-M,M]_\tau}\inf\limits_{s\in [t,t+L]_\tau}
e^{W(t)}}
{\eta \sum\limits_{k\in [-M,M]_\eta} e^{
W(\eta k)
}}},
$$
where big $M$ and small $\tau,\eta$ are appropriately chosen positive numbers.
In the following lemma we give a lower bound for $\mathcal{F}_{2H}(L)$.
\begin{lem}\label{parisian_pickands_constant_lower_bound}
It holds that for any $L>0$ and $H\in (0,1)$
$$\mathcal F_{2H}(L)
\ge  \E{
\left(\int\limits_\R e^{W(t)}dt\right)^{-1}} e^{-L^{2H}}
\sup\limits_{m>0}\left(e^{-\sqrt 2 m L^H}
\pk{\sup\limits_{s\in [0,1]}B_{H}(s)<m}\right)$$
with
$\E{
\left(\int\limits_\R e^{W(t)}dt\right)^{-1}}
 \in (0,\IF).$
\end{lem}
Taking $m = 1/\sqrt 2$ in the $\sup$ above we obtain a
useful for large $L$ estimate
$$\mathcal F_{2H}(L)\ge Ce^{-L^{2H}-L^H}, \ \ L>0$$
where $C$ is a some positive number that depends only on $H$.
The following lemma provides us an upper bound for the truncation error:

\begin{lem}\label{Pickands_constant_tail_estimate}
 For some fixed constant $c'>0$ and $M,L> 0$ it holds that
\bqny{
\left|\mathcal F_{2H}(L)-
\E{ \frac{\sup\limits_{t\in [-M,M]}\inf\limits_{s\in [t,t+L]}
e^{W(t)}}{\int\limits_{[-M,M]} e^{W(t)}dt}}
\right| \le e^{-c'M^{2H}}.
}
\end{lem}
Based on 2 lemmas above we propose the following algorithm for simulation
of $\mathcal F_{2H}(L)$:
\\\\
1) Take $M = \max(10L,5), \  \tau=\eta = 0.005$ and $n=10^4$;\\
2) simulate $n$ times $B_H(t), \ t \in [-M,M]_\tau$,
i.e, obtain $B_{H}^{(i)}(t), \ 1\le i\le n$;\\
3) calculate
 \bqny{
 \widehat{\mathcal F_{2H}(L)}:=  \frac{1}{n}\sum\limits_{i=1}^{n}
\frac{\sup\limits_{t\in [-M,M]_\tau}\inf\limits_{s\in [t,t+L]_\tau}
e^{\sqrt 2 B_{H}^{(i)}(s)-|s|^{2H}}}
{\eta \sum\limits_{k\in [-M,M]_\eta} e^{\sqrt 2 B_{H}
^{(i)}(k\eta)-|k\eta|^{2H}}}.
}

We give the proofs of all Lemmas above at the end of Section Proofs.

\section{Approximate values of Pickands \& Piterbarg constants}
In this section we apply
both algorithms introduced above and obtain
approximate numerical values for some particular choices
of parameters.
To implement our approach we use MATLAB software.
\\

\textbf{Piterbarg constant.}
We simulate several graphs of $\widehat{\mathcal F^h_L}$ for
different choices of $a$ and $b$:


\begin{figure}[H]
     \centering
     \begin{subfigure}[b]{0.3\textwidth}
         \centering
         \includegraphics[width=\textwidth]{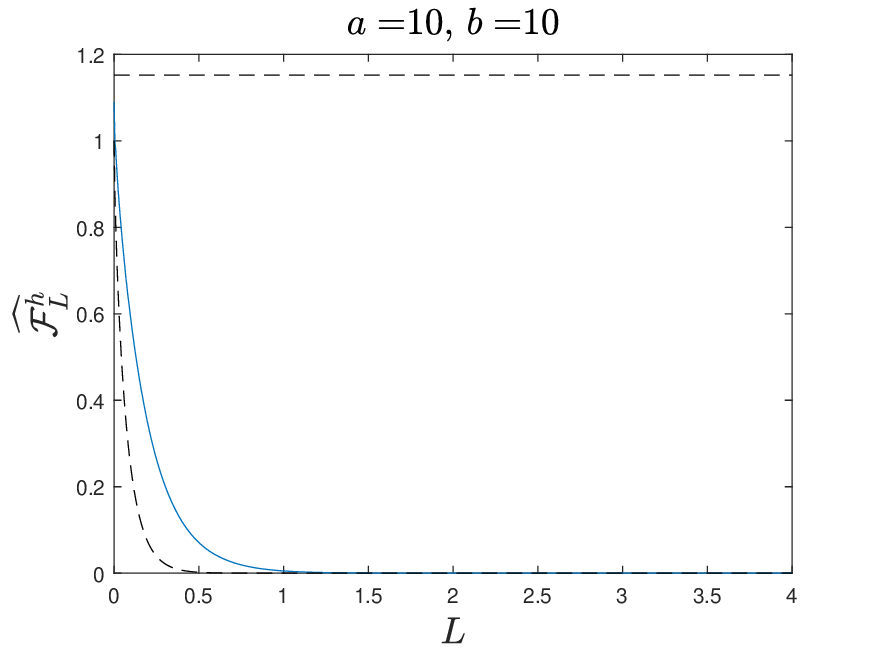}
     \end{subfigure}
     \hfill
     \begin{subfigure}[b]{0.3\textwidth}
         \centering
         \includegraphics[width=\textwidth]{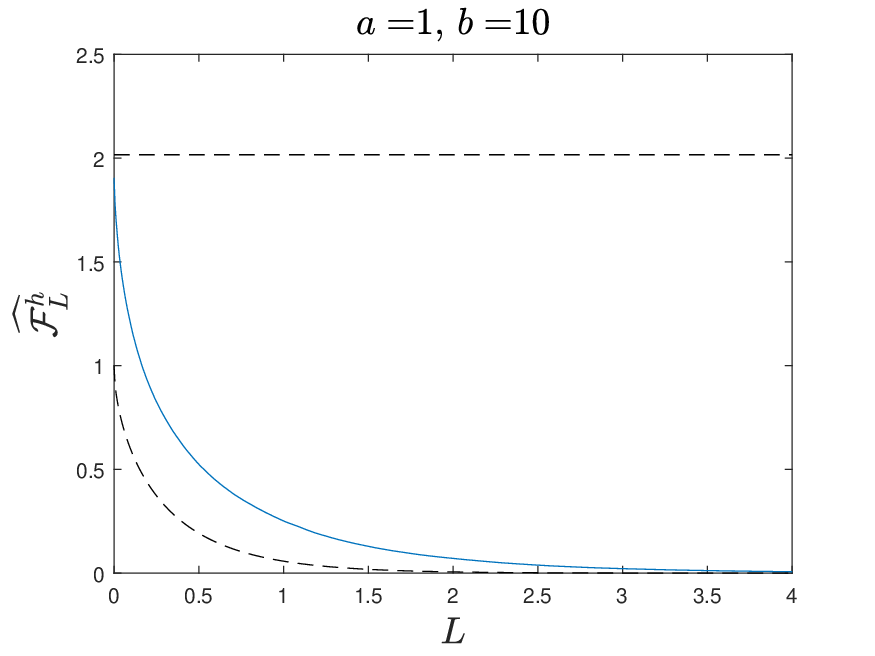}
     \end{subfigure}
     \hfill
     \begin{subfigure}[b]{0.3\textwidth}
         \centering
         \includegraphics[width=\textwidth]{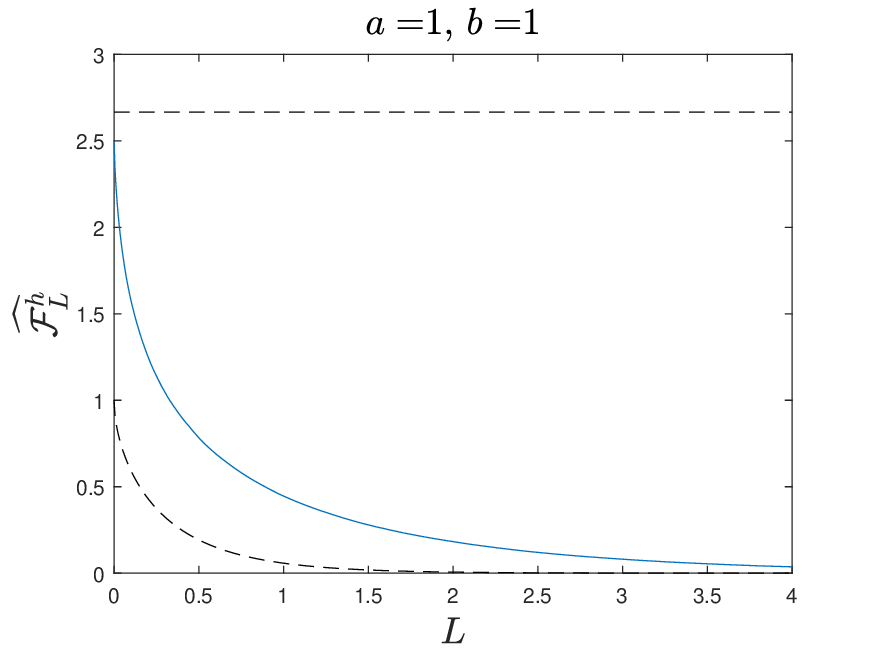}
     \end{subfigure}
\end{figure}

On each graph above the blue line is simulated value and
the dashed lines are theoretical bounds given in Lemma \ref{lower_bound_for_Piterb_const}.
We observe that the simulated values are between
the theoretical bounds given in Lemma \ref{lower_bound_for_Piterb_const},
$\widehat{\mathcal F_L^h}$ is decreasing function and
$\widehat{\mathcal F_L^h}$ tends to $1+\frac{1}{a}+
\frac{1}{b}-\frac{1}{a+b+1}$ as $L \to 0$.
\\

\textbf{Pickands constant.}
We plot several graphs of $\widehat{\mathcal F_{2H}(L)}$ for
different choices of $H$. Since the value of $\mathcal F_{1}(L)$
is known, namely,
$\mathcal F_{1}(L) =
\frac{e^{-L/4}-\sqrt{\pi L}\Phi(-\sqrt {L/2})}
{e^{-L/4}+\sqrt{\pi L}\Phi(\sqrt {L/2})}, \
L\ge 0, $
(see, e.g., \cite{LCPalowski})
we consider for simulations only cases $H<0.5$ and $H>0.5$.
To simulate fBm we use Choleski method, (see, e.g, \cite{MR1984656}).
\\





\emph{Short-range dependence case.}
We take $H$ equal to $0.1,0.2,0.3$ and $0.4$
and plot $\widehat{\mathcal F_{2H}(L)}$ for these values.

\begin{figure}[H]
\centering
	\includegraphics[scale = .7]{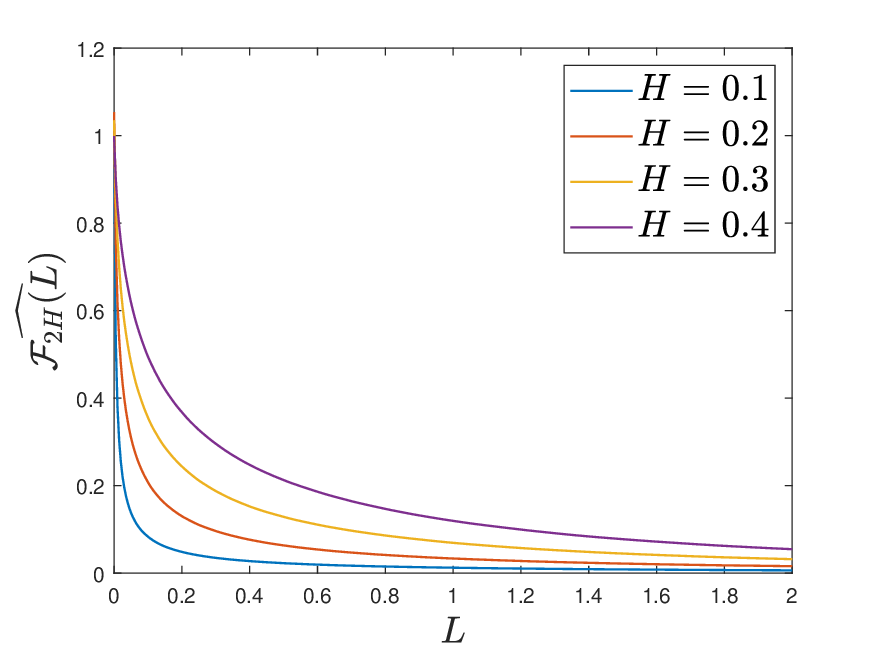}
\end{figure}

\emph{Long-range dependence case.}
Here we take  $H$ from $\{0.6,0.7,0.8,0.9\}$
and plot $\widehat{\mathcal F_{2H}(L)}$ for these values.
\begin{figure}[H]
\centering
	\includegraphics[scale = .7]{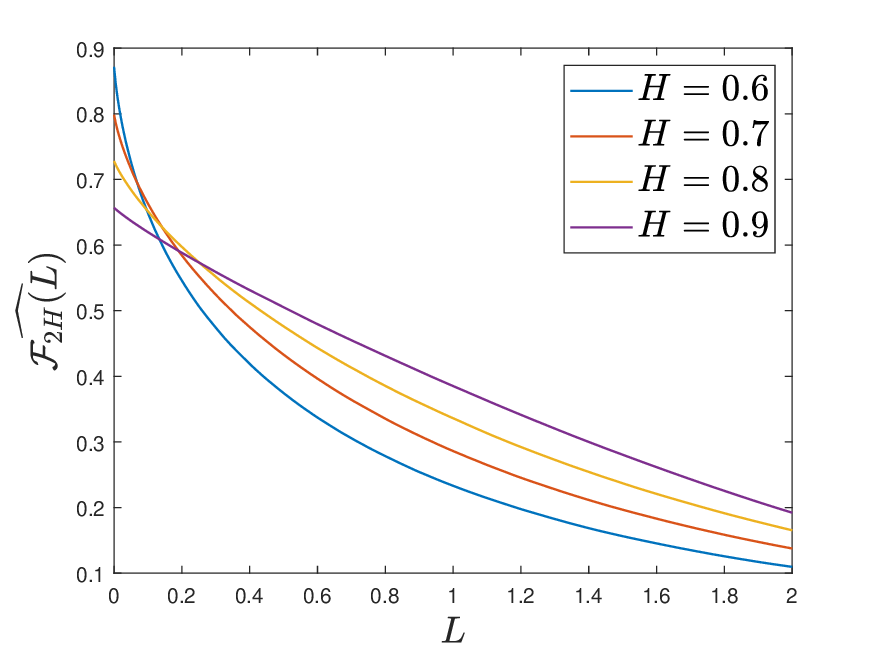}
\end{figure}

Observe that $\widehat{\mathcal F_{2H}(L)}$
is a strictly decreasing function of $L$ for all $H\in (0,1)$.
It seems also that $\mathcal F_{2H}(L)$ for fixed $L$ is an increasing
function of $H$ for $H\in (0,1/2)$ and is not an increasing function
of $H$ for $H\in (1/2,1)$.

\section{Proofs}
Before giving our proofs we formulate a few auxiliary statements.
As shown, e.g., in \cite{20lectures}
\bqn{\label{ratio}
\PH(x) \sim \frac{e^{-x^2/2}}{\sqrt {2\pi}x}, \ \ x \to \IF.
}

Recall that $K_H,D_1$ and $\mathbb{C}_H^{(1)}$ are defined in
\eqref{mathbb_D}.
The following result immediately follows from \cite{ParisRuinGenrealHashorvaJiDebicki,LCPalowski}:
\begin{prop}\label{theo_app}
Assume that $T_u$ satisfies \eqref{assumption_T_u}.
Then as $u\to \IF$
\bqny{
\pk{\sup\limits_{t \ge 0}\inf\limits_{[t,t+T_u]}(B_H(t)-c_1t)>q_1u} \sim
\begin{cases}
\frac{e^{-c_1^2T/2}-c_1\sqrt{2\pi T}\Phi(-c_1\sqrt T)}
{e^{-c_1^2T/2}+c_1\sqrt{2\pi T}\Phi(c_1\sqrt T)}e^{-2c_1q_1u},
&H=1/2,\\
K_H\mathcal{F}_{2H}(TD_1)(\mathbb{C}_H^{(1)}u^{1-H})^{\frac{1}{H}-1}
\PH(\mathbb{C}_H^{(1)}u^{1-H}),
& H\neq 1/2.
\end{cases}}
\end{prop}

Now we are ready to present our proofs.
\\

\textbf{Proof of Theorems \ref{theoparis_simple} and \ref{theoparis}.} Since Theorem
\ref{theoparis_simple} follows immediately from Theorem
\ref{theoparis} we prove Theorem
\ref{theoparis} only. \\

\textbf{Case (1).} Assume that $t_*<t_1$.
Let
\bqny{
\psi_i(T_u,u) = \pk{\sup\limits_{t \ge 0}\inf\limits_{[t,t+T_u]}
(B_H(t)-c_it)>q_iu}, \ \ i=1,2.
}
For $0<\ve<t_1-t_*$ by the self-similarity of fBm
we have
\bqny{
\psi_1(T_u,u) \ge \pr(u)
&\ge&
\pk{\exists t\in (t_1-\ve,t_1+\ve):\inf
\limits_{s\in [t,t+T_u/u]}V_1(t)> u^{1-H},
\inf
\limits_{s\in [t,t+T_u/u]}V_2(t)> u^{1-H}}
\\ &=&
\pk{\exists t\in (t_1-\ve,t_1+\ve):\inf
\limits_{s\in [t,t+T_u/u]}V_1(t)> u^{1-H}}
,}
where
\bqny{
V_i(t) = \frac{B_H(t)}{c_it+q_i}, \ \ \ i =1,2.
}
We have by Borel-TIS inequality, see \cite{20lectures}
(details are in the Appendix)
\bqn{\label{mainintpariscase1}
\psi_1(T_u,u) \sim
\pk{\exists t\in (t_1-\ve,t_1+\ve):\inf
\limits_{s\in [t,t+T_u/u]}V_1(t)> u^{1-H}}, \ \ \ u \to \IF
}
implying $\pr(u)\sim \psi_1(T_u,u)$ as $u \to \IF$.
The asymptotics of $\psi_1(T_u,u)$ is given
in Proposition \ref{theo_app}, thus the claim follows.
\\
\\
Assume that $t_*=t_1$.
We have
\bqny{
\pk{\exists t\in [t_1,\IF):
\inf
\limits_{s\in [t,t+\frac{T_u}{u}]}V_1(s)> u^{1-H}}
&\le& \pr(u)
\\&\le&
\pk{\exists t\in [t_1,\IF):
\inf
\limits_{s\in [t,t+\frac{T_u}{u}]}V_1(s)> u^{1-H}}
\\& \ & + \ \
\pk{\exists t\in [0,t_1]:
V_2(t)> u^{1-H}}
.}
From the proof of Theorem 3.1,
case (4) in \cite{Lanpeng2BM} it follows that
the second term in the last line above
is negligible comparing with the final asymptotics of
$\pr(u)$ given in \eqref{parisclaim11}, hence
$$
\pr(u) \sim
\pk{\exists t\in [t_1,\IF):
\inf
\limits_{s\in [t,t+\frac{T_u}{u}]}V_1(s)> u^{1-H}},
\quad u \to \IF.
$$
By the same arguments as in \eqref{mainintpariscase1} it follows that
for $\ve>0$ the last probability above
is equivalent with
\bqny{
\pk{\exists t\in [t_1,t_1+\ve]:\inf
\limits_{s\in [t,t+T_u/u]}V_1(s)> u^{1-H}}, \ \ \ u \to \IF
.}

Since
$\mathcal{F}_{1}(T) = \frac{e^{-T/4}-\sqrt {\pi T}\Phi(-\sqrt {T/2}) }
{e^{-T/4}+\sqrt{\pi T}\Phi(\sqrt {T/2})}, \ T\ge 0$
(see \cite{ParisRuinGenrealHashorvaJiDebicki})
applying Theorem 3.3 in \cite{ParisianFiniteHoryzon}
with parameters in the notation therein
\bqny{
\widetilde{\sigma} = \frac{t_1^H}{c_1t+q_1}, \ \
\beta_1 = 2, \ \ D = \frac{1}{2t_1^{2H}}, \ \ \alpha = 2H, \ \
A = \frac{q_1^{H-3}H^{H-1}(1-H)^{4-H}}{2c_1^{H-2}}
}
we obtain
\bqny{
\pk{\exists t\in [t_1,t_1+\ve]:\inf
\limits_{s\in [t,t+T_u/u]}V_1(s)> u^{1-H}} \sim \frac{1}{2}
K_H\mathcal{F}_{2H}(T D_1)(\mathbb{C}_H^{(1)}u^{1-H})^{\frac{1}{H}-1}
\PH(\mathbb{C}_H^{(1)}u^{1-H}), \ \ \ \ u \to \IF
}
 and the claim is established.
Case $t_*\ge t_2$ follows by the same arguments.
\\
\\
\textbf{Case (2).} Define
\bqn{\label{Z_H}
Z_H(t) = \frac{B_H(t)}{\max(c_1t+q_1,c_2t+q_2)}, \quad t\ge 0.
}
Similarly to the proof of \eqref{mainintpariscase1}
we have by Borell-TIS inequality for $\ve>0$
as $u \to \IF$
\bqny{
 \pr(u)\notag
 &=& \pk{\exists t\ge 0: \inf\limits_{s\in [t,t+T_u/u]} Z_H(t)>u^{1-H}}
 \\&\sim&
 \pk{\exists t\in ( t_*-\ve,
 t_*+\ve) : \inf\limits_{s\in [t,t+T_u/u]} Z_H(t)>u^{1-H}}
=: p(u),\ \ u \to \IF.
}

Assume that $H<1/2$.
By "the double-sum" approach, see the proofs of
Theorem 3.1, Case (3) $H<1/2$ in \cite{Lanpeng2BM} and
Theorem 3.3. case i) in \cite{ParisianFiniteHoryzon} we have as $u \to \IF$
\bqn{\label{star}
p(u)
 \sim
\pk{\exists t\in (t_*,
 t_*+\ve)\!: \!\!\!\!\!\!\inf\limits_{s\in [t,t+\frac{T_u}{u}]}\!
  \!\!V_1(t)\!>\!u^{1-H}} \!+\! \pk{\exists t\in (
 t_*-\ve,t_*)\!:\!\!\!\!\!\! \inf\limits_{s\in [t,t+\frac{T_u}{u}]}
 \!\!\! V_2(t)\!>\!u^{1-H}}
.}
To compute the asymptotics of each probability in the line above
we apply Theorem 3.3 in \cite{ParisianFiniteHoryzon}.
For the first probability we have in the notation therein
$$ \widetilde{\sigma} = \frac{t_*^H}{c_1t_*+q_1}, \quad
\beta_1 = 1,\quad D = \frac{1}{2t_*^{2H}},\quad
\alpha = 2H<1, \quad A = \frac{t_*^{H-1}|H(c_1t_*+q_1)-c_1t_*|}{(c_1t_*
+q_1)^2}$$
implying as $ u \to \IF$
\bqny{ \pk{\exists t\in (
 t_*,  t_*+\ve) : \inf\limits_{s\in [t,t+\frac{T_u}{u}]} V_1(t)>u^{1-H}}
 \sim
\mathcal{F}_{2H}
(\frac{(c_1t_*+q_1)^{ \frac{1}{H} }
}{2^{\frac{1}{2H}}t_*^2} T)
\frac{ t_*^H\mD_H^{\frac{1}{H}-1}
u^{(1-H)(\frac{1}{H}-2)}}{|H(c_1t_*+q_1)-c_1t_*| 2^{\frac{1}{2H}}}
\PH(\mD_H u^{1-H}).}
Applying again
Theorem 3.3 in \cite{ParisianFiniteHoryzon}
we obtain the asymptotics of the second summand and the claim follows by
\eqref{star}.
\\
\\
Assume that $H=1/2$. In order to compute the asymptotics of $p(u)$
applying Theorem 3.3
in \cite{ParisianFiniteHoryzon}
with parameters
\bqny{ \alpha = \beta_1=\beta_2= 1, \ A_{\pm} = \frac{q_1-c_1t_*}{q_1+c_1t_*}, \
A = \frac{q_2-c_2t_*}{q_2+c_2t_*}, \ \widetilde{\sigma} = \frac{\sqrt {t_*}}
{c_1t_*+q_1}, \ D = \frac{1}{2t_*}
  }
we obtain ($d(\cdot)$ and $T'$ are defined in \eqref{defT'd})
\bqny{
p(u) \sim
\mathcal{F}_{T'}^d \PH(\mathbb{D}_{1/2} \sqrt u), \  \ u \to \IF
.}
Assume that $H>1/2$. Applying Theorem 3.3
in \cite{ParisianFiniteHoryzon} with parameters
$\alpha=2H>1=\beta_1=\beta_2$ we complete the proof since
$$p(u) \sim \PH(\mD_H u^{1-H}),\quad u \to \IF.  \ \ \ \  \ \
\  \ \ \ \ \ \ \ \Box$$

\textbf{Proof of Proposition \ref{remark_H<1/2}.}
\emph{Lower bound.}
Take $\kappa = 1-3H$ and recall that $\alpha = \frac{T^{2H}}{2t_*^{2H}}$.
We have
\bqn{
\mathcal{P}_T(u) &\ge&
\pk{\forall t \in [t_*-T/u,t_*] V_2(t)>u^{1-H}
\text{ and }
V_2(t_*)>u^{1-H}+\alpha u^{\kappa}}
\notag\\&\ge&\bar C \label{3}
\pk{V_2(t_*)>u^{1-H}+\alpha u^{\kappa}}
\\&\sim& \bar{C} \PH(\mD_H u^{1-H})e^{-C_{1,\alpha} u^{1-H+\kappa}-
 C_{2,\alpha}u^{2\kappa}}, \quad u \to \IF
\notag
,}
where $\bar C$ is a fixed positive constant that does not depend on $u$ and
$C_{1,\alpha}$ and $C_{2,\alpha}$ are defined in \eqref{C_alpha_1,2}.
Thus, to prove the lower bound we need to show \eqref{3}.
Note that $\eqref{3}$ is the same as
\bqny{
\pk{\exists t \in [t_*-T/u,t_*] :V_2(t)\le u^{1-H}
\text{ and }
V_2(t_*)>u^{1-H}+\alpha u^{\kappa}}
\le \ve'\pk{V_2(t_*)>u^{1-H}+\alpha u^{\kappa}},
}
with some $\ve'\in (0,1)$.
The last line above is equivalent with
\bqny{
& \ &
\pk{\exists t \in [ut_*-T,ut_*]: B_H(t)-c_2t\le
q_2u \text{ and }
B_H(ut_*)-c_2ut_*>q_2u+b\alpha u^{\kappa+H}}
\\&\le& \ve'\pk{B_H(ut_*)-c_2ut_*>q_2u+b\alpha u^{\kappa+H}},}
where $b = c_2t_*+q_2$.
We have with $\varphi_u(x)$ the density of $B_H(ut_*)$ that
the left part of the inequality above does not exceed
\bqny{& \ &
\pk{\exists t \in [ut_*-T,ut_*]:
B_H(ut_*)-B_H(t)>b\alpha  u^{\kappa+H}
\text{ and } B_H(ut_*)>b u}
 \\&=&
\int\limits_{bu}^\IF\pk{\exists t \in [ut_*-T,ut_*]:
x-B_H(t)>b\alpha  u^{\kappa+H}
 |B_H(ut_*)=x}\varphi_u(x)dx
\\&\le&
\int\limits_{bu}^{bu+1}\pk{\exists t \in [ut_*-T,ut_*]:
x-B_H(t)>b\alpha  u^{\kappa+H}
 |B_H(ut_*)=x}\varphi_u(x)dx
+\int\limits_{bu+1}^\IF\varphi_u(x)dx
.}
We also have that
$$\pk{B_H(ut_*)-c_2ut_*>q_2u} = \int\limits_{bu}^\IF\varphi_u(x)dx
\ge  \int\limits_{bu}^{bu+1}\varphi_u(x)dx
.$$
By \eqref{ratio} we have that
 $\int\limits_{bu+1}^\IF\varphi_u(x)dx$ is negligible comparing
with the last integral above.
Thus, to prove \eqref{3} we need to show
\bqny{
\int\limits_{bu}^{bu+1}\pk{\exists t \in [ut_*-T,ut_*]:
x-B_H(t)>b\alpha  u^{\kappa+H}
 |B_H(ut_*)=x}\varphi_u(x)dx
 \le \ve'\int\limits_{bu}^{bu+1}\varphi_u(x)dx, \ \ \ \ \ u \to \IF
,}
that follows from the inequality
\bqn{\label{borr}
\sup\limits_{x \in [bu,bu+1]}
\pk{\exists t \in [ut_*-T,ut_*]:
x-B_H(t)>b\alpha  u^{\kappa+H}
 |B_H(ut_*)=x} \le \ve'', \ \ \   u \to \IF,
}
where $\ve''\in (0,1)$ is some number.
We show the line above in the Appendix, thus the lower bound
holds.
\\\\

\emph{Upper bound.}
We have by the self-similarity of fBm
$$\mathcal{P}_T(u) = \pk{\sup\limits_{t\ge 0}\inf\limits_{s\in
[t,t+T/u]}Z_H(s)>u^{1-H}},$$
where $Z_H$ is defined in \eqref{Z_H}. For $\ve>0$ by
Borell-TIS inequality with
$I(t_*) = (-u^{-\ve}+t_*,t_*+u^{-\ve})$ we have
$$
\pk{\sup\limits_{t \notin I(t_*)}\inf\limits_{s\in
[t,t+T/u]}Z_H(s)>u^{1-H}} \le
\pk{\sup\limits_{t \notin I(t_*)}Z_H(t)>u^{1-H}}
\le \PH\left(\mD_Hu^{1-H}\right)e^{-Cu^{2-2H-2\ve}}, \quad u \to \IF,
$$
that is asymptotically smaller than the lower bound in
\eqref{paris_ineq_Hsmall} for sufficiently small $\ve$.
Thus, we focus on estimation of
$$q(u):=\pk{\sup\limits_{t \in I(t_*)}
\inf\limits_{s\in [t,t+T/u]}Z_H(s)>u^{1-H}}.$$
Denote $z^2(t) = \Var\{Z_H(t)\}$ and $\overline Z_H(t) = Z_H(t)/z(t)$.
By Lemma 2.3 in \cite{PickandsA} we have with $M = \max(z(t),z(t+T/u))$
(note, $1/M \ge \mD_H$)
\bqn{\label{mestim}
q(u)
&\le&  \pk{\exists t \in I(t_*): Z_H(t)>u^{1-H},
 Z_H(t+T/u)>u^{1-H} }
\notag\\&=&
\pk{\exists t \in I(t_*): \overline Z_H(t)>u^{1-H}/z(t),
\overline Z_H(t+T/u)>u^{1-H}/z(t+T/u) }
\notag\\&\le&
\pk{\exists t \in I(t_*): \overline Z_H(t)>u^{1-H}/M,
\overline Z_H(t+T/u)>u^{1-H}/M }
\notag\\&\le&2(1+o(1))
\PH\left(\frac{u^{1-H}}{M}\right)
\PH\left(\frac{u^{1-H}}{M} \sqrt{\frac{1-r(t,t+T/u)}{1+r(t,t+T/u)}}\right)
\notag\\&\le&2(1+o(1))
\PH\left(\frac{ u^{1-H}}{M}\right)
\PH\left(\mD_H u^{1-H} \sqrt{\frac{1-r(t,t+T/u)}{2}}\right),
}
where $r$ is the correlation function of $Z_H$.
Since $r(t,s) = corr (B_H(t),B_H(s))$ we have for all $t\in I(t_*)$
\bqny{
1-r(t,t+T/u) = \frac{T^{2H}}{2t_*^{2H}}u^{-2H}+O\left(u^{-2H}(|t-t_*|
+|t+T/u-t_*|)+u^{-2}\right)
,\quad u \to \IF}
implying
\bqny{
\mD_H u^{1-H} \sqrt{\frac{1-r(t,t+T/u)}{2}} =
u^{1-2H}\frac{T^H\mD_H}{2t_*^H}+O(u^{1-2H}(|t-t_*|+|t+T/u-t_*|)+u^{-1})
,\quad u \to \IF.}
Thus, by \eqref{ratio} we obtain
\bqn{\label{asym1}
\PH\left(\mD_H u^{1-H} \sqrt{\frac{1-r(t,t+T/u)}{2}}\right) \le
 \PH\left(u^{1-2H}\frac{T^H\mD_H}{2t_*^H}\right)
 e^{Cu^{2-4H}(|t-t_*|+|t+T/u-t_*|)}
 ,\quad u \to \IF.
}
Next we have as $u \to \IF$ for some $C_1>0$
\bqny{
\PH\left(\frac{ u^{1-H}}{M}\right) \sim \PH(\mD_H u^{1-H})e^{-C_1u^{2-2H}(|t-t_*|+|t+T/u-t_*|)}
}
and by \eqref{asym1} we have for all $t\in I(t_*)$ and large $u$
$$\PH\left(\frac{ u^{1-H}}{M}\right)\PH\left(\mD_H u^{1-H} \sqrt{\frac{1-r(t,t+T/u)}{2}}\right)
\le
\PH\left(\mD_H u^{1-H}\right) \PH\left(u^{1-2H}\frac{T^H\mD_H}{2t_*^H}\right)
 e^{(Cu^{2-4H}-C_1u^{2-2H})(|t-t_*|+|t+T/u-t_*|)}
 $$
and the claim follows from the line above and
\eqref{mestim}.
\QED\\\\

\textbf{Proof of Lemma \ref{lower_bound_for_Piterb_const}.}
\emph{Lower bound.} We have
\bqny{
\sup\limits_{t\in \R}\inf\limits_{s\in [t,t+L]}
e^{\sqrt 2 B(s)-|s|+h(s)} \ge
\inf\limits_{s\in [0,L]}
e^{\sqrt 2 B(s)-(1+a)s}
\ge
e^{-(1+a)L}\inf\limits_{s\in [0,L]}e^{\sqrt 2 B(s)}
\overset{d}{=} e^{-(1+a)L} e^{-\sup\limits_{s\in [0,L]}\sqrt 2 B(s)},
}
where the symbol '$\overset{d}{=}$' means equality in
distribution between two random variables.
Taking expectations of both sides in the line above we obtain
\bqny{
\mathcal{F}_L^{h} \ge e^{- L(1+a)}
\E{e^{-\sup\limits_{s\in [0,L]}\sqrt 2 B(s)}},
}
and our next step is to calculate the expectation above.
It is known (see, e.g., Chapter 11.1 in \cite{20lectures})
that
$$\pk{\sup\limits_{s\in [0,L]}\sqrt 2 B(s)>x}
= 2\pk{\sqrt 2 B(L)>x} = 2\PH\left(\frac{x}{\sqrt{2L}}\right),\ x>0$$
hence we obtain that
$\frac{e^{-x^2/4L}}{\sqrt{\pi L}}, \  x>0$ is
the density of $\sup\limits_{s\in [0,L]}\sqrt 2 B(s)$.
Thus, we have
\bqny{
\E{e^{-\sup\limits_{s\in [0,L]}\sqrt 2 B(s)}} =
\int\limits_0^\IF e^{-x}\frac{e^{-x^2/4L}}{\sqrt{\pi L}}dx =
\frac{e^{L}}{\sqrt{\pi L}}\int\limits_0^\IF
e^{-(\frac{x}{2\sqrt L}+\sqrt L)^2}dx =
\frac{2 e^{L}}{\sqrt{\pi }}\int\limits_{\sqrt L}^\IF
e^{-z^2}dz
 = 2e^{L}\PH(\sqrt {2L}),
}
and combining all calculations above we obtain
$$\mathcal{F}_L^{h} \ge 2e^{- La}\PH(\sqrt{2L}), \ \ \ L\ge 0.$$
On the other hand we have
\bqny{
\sup\limits_{t\in \R}\inf\limits_{s\in [t,t+L]}
e^{\sqrt 2 B(s)-|s|+h(s)} \ge \inf\limits_{s\in [-L,0]}
e^{\sqrt 2 B(s)-(1+b)|s|} \overset{d}{=} \inf\limits_{s\in [0,L]}
e^{\sqrt 2 B(s)-(1+b)s},
}
and estimating $\inf \limits_{s\in [0,L]}
e^{\sqrt 2 B(s)-(1+b)s}$ as above we have
$\mathcal{F}_L^{h} \ge 2e^{- Lb}\PH(\sqrt{2L}), \ L\ge 0$, that completes
 the proof of the lower bound.
\\\\
\emph{Upper bound.} Note that $\mathcal{F}_{2H}^L\le \mathcal F_{2H}^0$ and
hence since a BM has independent branches for positive and
negative time we have with $B_*$ an independent BM
\bqny{
\mathcal{F}_{2H}^L\le
\E{\sup\limits_{t\inr} e^{\sqrt 2 B(t)-h(t)} }
&=&
\E{\max\left(\sup\limits_{t\ge 0} e^{\sqrt 2 B(t)-(a+1)t}, \,
\sup\limits_{t\le 0} e^{\sqrt 2 B(t)-(b+1)|t|}
\right)}
\\&=&
\E{\max\left(\sup\limits_{t\ge 0} e^{\sqrt 2 B(t)-(a+1)t}, \,
\sup\limits_{t\ge 0} e^{\sqrt 2 B^*(t)-(b+1)t}
\right)} = \E{e^{\max(\xi_a,\xi_b)}},
}
where $\xi_a$ and $\xi_b$ are exponential random variables with survival
functions $e^{-(a+1)x}$ and $ e^{-(b+1)x}$, respectively,
see \cite{DeM15}. Since $\xi_a$ and
$\xi_b$ have exponential distributions the last expectation above is
$1+\frac{1}{a}+\frac{1}{b}-\frac{1}{a+b+1}$
and the claim follows.
\QED \\\\

\textbf{Proof of Lemma \ref{lemma_estimation_of_supinf_over_infty}.}
First we have
\bqny{
\E{\sup\limits_{t\in \R \backslash [-M,M]}\inf\limits_{s\in [t,t+L]}
e^{\sqrt 2 B(s)-|s|+h(s)}}
 \le
\E{\sup\limits_{s\in [M,\IF)}
e^{\sqrt 2 B(s)-(a+1)s}}
+ \E{\sup\limits_{s\in (-\IF,-M]}e^{\sqrt 2 B(s)-(b+1)|s|}}.
}

Later on we work with the first expectation above. We have
\bqny{
& \ &
\E{\sup\limits_{s\in [M,\IF)} e^{\sqrt 2 B(s)-(1+a)s}}
\\&=&
\int\limits_\R e^x
\pk{\sup\limits_{s\in [M,\IF)} (\sqrt 2 B(s)-(1+a)s)>x}dx
\\&=&
\int\limits_\R e^x
\pk{\sup\limits_{s\in [M,\IF)} (\sqrt 2 (B(s)-B(M))-(1+a)(s-M))
>x+M(1+a)-\sqrt 2 B(M)}dx
.}
Since a BM has
independent increments we have with $B^*$ an independent BM that the
last integral above equals
\bqny{& \ &
\int\limits_\R e^x
\pk{\sup\limits_{s\in [0,\IF)} (\sqrt 2 B(s)-(1+a)s)
>x+M(1+a)-\sqrt {2M} B^*(1)}dx
\\&=&\frac{1}{\sqrt{2\pi}}\int\limits_\R
\int\limits_\R e^x e^{-z^2/2}
\pk{\sup\limits_{s\in [0,\IF)} (\sqrt 2 B(s)-(1+a)s)
>x+M(1+a)-\sqrt {2M} z}dx dz
.}
We know that $\pk{\sup\limits_{t\ge 0}(B(t)-ct)>x} = \min (1, e^{-2cx})$
for $c>0$ and $x\inr$, thus the expression above equals
\bqny{& \ &
\frac{1}{\sqrt{2\pi}}\int\limits_\R
\int\limits_\R e^{x-z^2/2}
\min(1, e^{-(1+a)(x+M(1+a)-\sqrt{2M}z)})dx dz
\\&=&\frac{1}{\sqrt{2\pi}}\int\limits_\R
\int\limits_{\frac{(1+a)M+x}{\sqrt{2M}}}^{\IF} e^{x-z^2/2}
dz dx+
\frac{1}{\sqrt{2\pi}}\int\limits_\R
\int\limits_{-\IF}^{\frac{(1+a)M+x}{\sqrt{2M}}}
 e^{x-z^2/2-(1+a)(x+M(1+a)-\sqrt{2M}z)}
dz dx
\\&=&\int\limits_\R
 e^{x}\PH(\frac{(1+a)M+x}{\sqrt{2M}}) dx+
\frac{1}{\sqrt{2\pi}}\int\limits_\R e^{-ax}
\int\limits_{-\IF}^{\frac{(1+a)M+x}{\sqrt{2M}}}
e^{-\frac{(z-\sqrt {2M}(1+a))^2}{2}}dzdx
\\&=&\int\limits_\R
 e^{x}\PH(\frac{(1+a)M+x}{\sqrt{2M}}) dx+
\frac{1}{\sqrt{2\pi}}\int\limits_\R e^{-ax}
\int\limits_{-\IF}^{\frac{-(1+a)M+x}{\sqrt{2M}}}
e^{-\frac{z^2}{2}}dzdx
\\&=&\int\limits_\R
 e^{x}\PH(\frac{(1+a)M+x}{\sqrt{2M}}) dx+
\int\limits_\R e^{-ax}
\Phi(\frac{-(1+a)M+x}{\sqrt{2M}})
dx.
}

Integrating the first integral above by parts we have
\bqny{
\int\limits_\R
 e^{x}\PH(\frac{(1+a)M+x}{\sqrt{2M}}) dx
= -\int\limits_\R
 \big(\PH(\frac{(1+a)M+x}{\sqrt{2M}})\big)'e^xdx
&=&
\frac{1}{\sqrt{2\pi}\sqrt{2M}}\int\limits_\R e^{-\frac{((1+a)M+x)^2}{4M}}
e^xdx
\\&=&
\frac{e^{-aM}}{\sqrt{2\pi}\sqrt{2M}}\int\limits_\R
e^{-\frac{((a-1)M+x)^2}{4M}}dx
=  e^{-aM}.
}
For the second integral we have similarly
\bqny{
\int\limits_\R e^{-ax}\Phi(\frac{-(1+a)M+x}{\sqrt{2M}})dx
&=&
-\frac{1}{a}\int\limits_\R \Phi\big(\frac{-(1+a)M+x}{\sqrt{2M}}\big)'
e^{-ax}dx
\\&=&
\frac{1}{a}\frac{1}{\sqrt{2\pi}\sqrt{2M}}
\int\limits_\R e^{-\frac{(-(1+a)M+x)^2}{4M}-ax}dx
=
\frac{e^{-aM}}{a\sqrt{2\pi}\sqrt{2M}}
\int\limits_\R e^{-\frac{((1-a)M+x)^2}{4M}}dx
=\frac{e^{-aM}}{a}
.}
Summarizing all calculations above we obtain
\bqny{
\E{\sup\limits_{t\in [M,\IF)}
e^{\sqrt 2 B(t)-(1+a)t}} = e^{-aM}\left(1+\frac{1}{a}\right).
}
By the same approach and the symmetry of BM around zero we have
\bqny{
\E{\sup\limits_{t\in (-\IF,-M]}
e^{\sqrt 2 B(t)-(1+b)|t|}} = e^{-bM}\left(1+\frac{1}{b}\right)
}
and hence combining both equations above with
the first inequality in the proof we obtain the claim. \QED
\\

\textbf{Proof of Lemma \ref{parisian_pickands_constant_lower_bound}.}
From \cite{DiekerY} it follows, that for any $L\ge 0$
\bqn{\label{pick_constant_formula} \mathcal F_{2H}(L) =
\E{ \frac{\sup\limits_{t\in \R}\inf\limits_{s\in [t,t+L]}
e^{W(s)}}{\int\limits_\R e^{W(t)}dt}}
,}
later on we use this formula in the proof. Observe that
$\sup\limits_{t\inr}\inf\limits_{s\in [t,t+L]}
e^{W(s)}
\ge \inf\limits_{s\in [0,L]}e^{W(s)}$,
 hence
\bqny{
\mathcal F_{2H}(L) \ge
\E{ \frac{\inf\limits_{s\in [0,L]}e^{W(s)
}}{\int\limits_\R e^{W(t)}dt}}
\ge
e^{-L^{2H}}\E{ \frac{e^{-\sqrt 2\sup\limits_{s\in [0,L]}B_{H}(s)}}
{\int\limits_\R e^{W(t)}dt}}
.}
Let $\xi = \sup\limits_{s\in [0,L]}B_{H}(s)$,
$(\Omega,\mathbb P)$ be the general probability space and
$\Omega_m = \{\omega\in \Omega: \xi(\omega)<m\}$ for $m>0$.
The last expectation above equals
\bqny{
\E{
\frac{e^{-\sqrt 2\xi}}{\int\limits_\R e^{W(t)}dt}}
&=& \int\limits_\Omega \frac{e^{-\sqrt 2 \xi(\omega)}}
{\int\limits_\R e^{\sqrt 2 B_{H}(t,\omega)-|t|^{2H}}dt} d\mathbb P(\omega)
\\&\ge&
\int\limits_{\Omega_m} \frac{e^{-\sqrt 2 \xi(\omega)}}
{\int\limits_\R e^{\sqrt 2 B_{H}(t,\omega)-|t|^{2H}}dt} d\mathbb P(\omega)
\\&\ge& \pk{\Omega_m}e^{-\sqrt 2 m }
\int\limits_{\Omega_m} \frac{1}
{\int\limits_\R e^{W(t)}dt}d\mathbb P(\omega)
\\&\ge&
e^{-\sqrt 2 m}\pk{\xi<m}
\E{\frac{1}{\int\limits_\R e^{W(t)}dt}}
.}
Next taking $m = n L^H$  by the self-similarity of fBm we have that
\bqny{
e^{-\sqrt 2 m}\pk{\xi<m} = e^{-\sqrt 2 nL^H }
\pk{\sup\limits_{s\in [0,L]}B_{H}(s)<nL^H}=
e^{-\sqrt 2 nL^H }
\pk{\sup\limits_{s\in [0,1]}B_{H}(s)<n}
.}
Taking $\sup$ with respect to $n$ over $(0,\IF)$ we have
$$\mathcal F_{2H}(L) \ge
\E{ \frac{1}
{\int\limits_\R e^{W(t)}dt}} e^{-L^{2H}}
\sup\limits_{n>0}\Big(e^{-\sqrt 2 n L^H}
\pk{\sup\limits_{s\in [0,1]}B_{H}(s)<n}\Big)
$$
and hence to complete the proof we need to show that
the expectation in the expression above is a finite
positive constant. Since the classical Pickands constant is finite
(see, e.g., \cite{PickandsA,DI2005,20lectures,DiekerY}) we have
\bqny{
0<\E{ \frac{1} {\int\limits_\R e^{W(t)}dt}}
\le
\E{ \frac{\sup\limits_{t\inr} e^{W(t)}}
{\int\limits_\R e^{W(t)}dt}} =\mathbb{H}_{2H} \in (0,\IF).
 \ \ \ \  \ \ \ \  \ \ \ \ \  \ \Box
}

\textbf{Proof of Lemma \ref{Pickands_constant_tail_estimate}.}
By \eqref{pick_constant_formula} we have that
\bqny{& \ &
\Big|\mathcal F_{2H}(L)-
\E{ \frac{\sup\limits_{t\in [-M,M]}\inf\limits_{s\in [t,t+L]}
e^{W(s)}}{\int\limits_{[-M,M]}
e^{W(t)}dt}}
\Big|
\\&=&\Big|\Big(
\E{ \frac{\sup\limits_{t\inr}\inf\limits_{s\in [t,t+L]}e^{W(s)}}{\int\limits_{\R} e^{W(t)}dt}}
-
\E{ \frac{\sup\limits_{t\in [-M,M]}\inf\limits_{s\in [t,t+L]}
e^{W(s)}}{\int\limits_{\R} e^{W(t)}dt}}\Big)
\\&\ &+\Big(
\E{ \frac{\sup\limits_{t\in [-M,M]}\inf\limits_{s\in [t,t+L]}
e^{W(s)}}{\int\limits_{\R} e^{W(t)}dt}}
-\E{ \frac{\sup\limits_{t\in [-M,M]}\inf\limits_{s\in [t,t+L]}e^{W(s)}}{\int\limits_{[-M,M]} e^{W(t)}dt}}
\Big)\Big|
\\&\le&
\E{ \frac{\sup\limits_{t\in \R\backslash[-M,M]}e^{W(t)}}{\int\limits_{\R}
e^{W(t)}dt}}
+
\E{
\sup\limits_{t\in [-M,M]}e^{W(t)}\frac{
\int\limits_{\R\backslash [-M,M]} e^{W(t)}dt}
{\int\limits_{\R} e^{W(t)}dt
\int\limits_{[-M,M]} e^{W(t)}dt
}}.}
As follows from Section 4 in \cite{DiekerY},
the last line above does not exceed $e^{-c'M^{2H}},$ and the claim holds.
\QED\\\\

\section{Appendix}

\textbf{Proof of \eqref{mainintpariscase1}.}
To establish the claim we need to show, that
\bqny{
\pk{\exists t\in \R \backslash[t_1-\ve,t_1+\ve]:\inf
\limits_{s\in [t,t+T/u]}V_1(s)> u^{1-H}} = o(\psi_1(T_u,u)), \ \
u \to \IF.
}
Applying Borell-TIS inequality (see, e.g., \cite{20lectures}) we have
as $u \to \IF$
\bqny{
\pk{\exists t\in \R \backslash[t_1-\ve,t_1+\ve]:\inf
\limits_{s\in [t,t+T/u]}V_1(s)> u^{1-H}} \le
\pk{\exists t\in \R \backslash[t_1-\ve,t_1+\ve]:
V_1(t)> u^{1-H}}
\le
e^{-\frac{(u^{1-H}-\widetilde M)^2}{2m^2}},
}
where
$$\widetilde M = \E{\sup\limits_{\exists t\in \R \backslash[t_1-\ve,t_1+\ve]} V_1(t)}
< \IF, \quad
m^2 = \max\limits_{\exists t\in \R \backslash[t_1-\ve,t_1+\ve]
 }\Var\{V_1(t)\}.$$
Since $\Var\{V_1(t)\}$ achieves its unique maxima at $t_1$ we obtain
by \eqref{ratio}
that
 $$e^{-\frac{(u^{1-H}-\widetilde M)^2}{2m^2}}  = o(\pk{V_1(t_1)<u^{1-H}}),
\quad u \to \IF$$
and the claim follows from the asymptotics of $\psi_1(T_u,u)$
given in Proposition \ref{theo_app}.
\QED
\\
\\
\textbf{Proof of \eqref{borr}.}
Define $X_{x,u}(t)= x-B_H(t)|B_H(ut_*) = x, \ t \in [ut_*-T,u]$.
To calculate the covariance and expectation of $X_{x,u}$ we use the formulas
$$\cov((B,C)|A = x)  = \cov(B,C) - \frac{\cov(A,B)\cov(A,C)}{\Var \{A\}}
\quad \text{and} \quad \E{B|A=x} = x\cdot \frac{\cov(A,B)}{\Var\{A\} },$$
where $A,B$ and $C$ are centered Gaussian random variables and
$x\inr$. We have for $x\in [bu,bu+1]$ and $ t,s\in [ut_*-T,ut_*]$
with $v=ut_*$, $y = 1-\frac{t}{v}$ and $z = 1-\frac{s}{v}$ as $u \to \IF$
\bqn{\label{cov}\notag & \ &
\cov(X_{x,u}(t),X_{x,u}(s))
\\&= &
\frac{t^{2H}+s^{2H}-|t-s|^{2H}}{2}-\frac{(t^{2H}+v^{2H}-|t-v|^{2H})
(s^{2H}+v^{2H}-|s-v|^{2H})}{4v^{2H}}
\notag\\&=&
\frac{v^{2H}}{4}\Big( 2(\frac{t}{v})^{2H}+2(\frac{s}{v})^{2H}-2
|\frac{t}{v}-\frac{s}{v}|^{2H}
-((\frac{t}{v})^{2H}+1-|\frac{t}{v}-1|^{2H})
((\frac{s}{v})^{2H}+1-|\frac{s}{v}-1|^{2H})
\Big)
\notag\\&=&
\frac{v^{2H}}{4}\Big(
2(1-y)^{2H}+2(1-z)^{2H}-2
|y-z|^{2H}-((1-y)^{2H}+1-y^{2H})
((1-z)^{2H}+1-z^{2H})\Big)
\notag\\&=&
\frac{v^{2H}}{4}\Big(
2-4Hy+2-4Hz+O(y^2+z^2)-2|y-z|^{2H}
\notag\\& \ & \quad \ \ \ \ - \
(2-2Hy-y^{2H}+O(y^2))(2-2Hz-z^{2H}+O(z^2))\Big)
\notag\\&=&
\frac{v^{2H}}{4}\Big(
2y^{2H}+2z^{2H}-2|y-z|^{2H}+O(y^2+z^2+z^{2H}y^{2H})
\Big)
\notag\\&=&
(1+o(1))\frac{(ut_*-t)^{2H}+(ut_*-s)^{2H}-|t-s|^{2H}}{2}
.}
For the expectation we have as $u \to \IF$
\bqny{
\E{X_{x,u}(t)} = x(1-\frac{v^{2H}+t^{2H}-|v-t|^{2H}}{2v^{2H}})
&=&\frac{x}{2}(1-(t/v)^{2H}+(1-t/v)^{2H})
\notag\\&\le& \frac{1}{2}(bu+1)(1 - (1-y)^{2H} + y^{2H})
\notag\\&\le& (bu/2+1)(1- 1+2H y -o(y) +y^{2H})
\notag\\&\le& Hbuy+\frac{1}{2}buy^{2H}+o(1)
.}
From the line above it follows that for some $C_*>0,$
$H<1/2, \ x\in [bu,bu+1]$ and $ t\in [ut_*-T,ut_*]$
\bqny{
\E{X_{x,u}(t)} \le C_*+\frac{u^{1-2H}b}{2t_*^{2H}}(ut_*-t)^{2H}.
}
We have
\bqny{& \ &
\sup\limits_{x\in[bu,bu+1]}\pk{\exists t \in [ut_*-T,ut_*]:X_{x,u}(t)
>u^{H+\kappa}\alpha b}
\\&=&
\sup\limits_{x\in[bu,bu+1]}
\pk{\exists t \in [ut_*-T,ut_*]:X_{x,u}(t)-\E{X_{x,u}(t)}
>u^{H+\kappa}\alpha b-\E{X_{x,u}(t)}}
\\&\le&
\pk{\exists t \in [0,T]: Y_u(t)+ f(t)>0},
}
where $Y_u(t) = X_{x,u}(ut_*-T+t)-\E{X_{x,u}(ut_*-T+t)},\ t\in[0,T]$ and
$f(t)$ is the linear
function such that $ f(T) = C_1$ and $ f(0) = -C_*<0$.
Next we have by \eqref{cov} for all large $u$ and $t,s \in
[0,T]$
\bqny{ & \ &
\E{(Y_u(t)+ f(t)-Y_u(s)- f(s))^2}
\\&=&
\E{(Y_u(t)-Y_u(s))^2}+C(t-s)^2
\\&\le&
C_1\Big( (ut_*-t)^{2H}+(ut_*-s)^{2H}-
(ut_*-t)^{2H}-(ut_*-s)^{2H}+|t-s|^{2H} \Big) +C(t-s)^2
\\&\le& 2|t-s|^{2H}.
}
Thus, by Proposition 9.2.4 in \cite{20lectures}
the family $Y_u(t)+ f(t), \ u>0,  \ t \in [0,T]$ is tight in
$\mathcal{B}(C([0,T]))$.
As follows from \eqref{cov}, it holds that
$\{Y_u(t)+ f(t)\}_{t \in [0,T]}$
converges to $\{B_H(t)+f(t)\}_{t \in [0,T]}$ in
the sense of convergence of finite-dimensional distributions as $u \to \IF$.
Hence by Theorems 4 and 5 in Chapter 5 in \cite{BylinskiiShiryaevBook}
the tightness and convergence of finite-dimensional distributions imply
weak convergence
$$\{Y_u(t)+ f(t)\}_{ t \in [0,T]}
 \Rightarrow \{B(t)+f(t)\}_{t \in [0,T]}.$$
Since the functional
$F(g) = \sup\limits_{t \in [0,T]}g(t)$
is continuous in the uniform metric we obtain
$$\pk{\exists t \in [0,T]: Y_u(t)+ f(t)>0}
\to
\pk{\exists t \in [0,T]: B_H(t)+ f(t)>0}, \ \ u \to \IF.
$$
Thus, to prove the claim it is enough to show that
\bqn{\label{fbm}
\pk{\exists t \in [0,T]: B_H(t)+ f(t)>0}<1.
}
We have for some large $m$ with $l(s)$ the density of $B_H(T)$
\bqn{\label{in}\notag
\pk{\sup\limits_{t\in[0,T]} (B_H(t)+ f(t))<0} &\ge&
\pk{\sup\limits_{t\in[0,T]} (B_H(t)+ f(t))<0 \text{ and }
B_H(T)<-m}
 \\&=&
\int\limits_{-\infty}^{-m}\pk{\sup\limits_{t\in[0,T]} (B_H(t)+ f(t))<0 |
B_H(T)=s}l(s)ds.
}
Define process
$\widetilde{B}_s(t) = B_H(t)+f(t)|B_H(T) = s,\ t\in[0,T]$.
We have for  $s<-m$ and $t\in [0,T]$
$$\E{\widetilde{B}_s(t)} = f(t)+s\frac{t^{2H}+T^{2H}-|T-t|^{2H}}{2T^{2H}}<
-C_1/2
,\quad \Var\{\widetilde{B}_s(t)\} = t^{2H}-
\frac{(T^{2H}+t^{2H}-|t-s|^{2H})^2}{4T^{2H}}<C_2$$
and thus
\bqny{
\pk{\sup\limits_{t\in[0,T]} (B_H(t)+ f(t))<0 |B_H(T)=s}
\ge
\pk{\sup\limits_{t\in[0,T]}\big(
\widetilde{B}_s(t)-\E{\widetilde{B}_s(t)}\big)<C_1/2}.
}
The last probability above is positive for any $s<-m$,
see Chapters 10 and 11 in \cite{LifBookLan2016} and hence the integral in
\eqref{in} is positive implying
$$\pk{\sup\limits_{t\in[0,T]} (B_H(t)+ f(t))<0}>0.$$
Consequently \eqref{fbm} holds and the claim
is established.
\QED
\\
\

\COM{
\textbf{Choice of $M$ for simulation of Piterbarg
constant.} We consider 2 cases: $L<\frac{1}{4\pi}$ and $L\ge \frac{1}{4\pi}$. \\
\emph{Assume that $L<\frac{1}{4\pi}$.} Then $
\PH(\sqrt{2L})\in (\PH(\frac{1}{\sqrt{2\pi}}),1/2)$ hence we can neglect it. We need to have
$$\min(a,b)M-L\min(a,b)>z$$ in order
to have maximal possible percentage of the error $\approx e^{-z}$.
Taking $z=7$ we can make the error smaller
then $10^{-3}$; we do not need to have better accuracy since there is
also the error of discretization. Thus, inequality above implies choice of $M$ as
$$ M = \frac{7+L(2+\min(a,b))}{\min(a,b)} .$$
\emph{Assume that $L\ge \frac{1}{4\pi}$. } In this case we have
$\PH(\sqrt{2L}) \le \frac{1}{2\sqrt{\pi L}}e^{-L}\le e^{-L}$ and thus
we need to have
$e^{-\min(a,b)M}<< e^{-L(3+\min(a,b))}$
and our choice is
$$M := 1+\frac{7+3L}{\min(a,b)}.$$
Since if $L<\frac{1}{4\pi}$, then
$\frac{7+L(2+\min(a,b))}{\min(a,b)} \approx
\frac{7+L(3+\min(a,b))}{\min(a,b)}$, and in the following
we take $M  = \frac{7+L(3+\min(a,b))}{\min(a,b)}$.\\
}

{\bf Acknowledgement:}
Grigori Jasnovidov was supported by Ministry of Science and Higher Education of the Russian Federation grant 075-15-2022-289.

\bibliography{EEEA}{}
\bibliographystyle{unsrt}
\end{document}